\documentclass[a4paper,10pt]{amsart}
\usepackage[arrow,matrix]{xy}
\usepackage{amsmath,amssymb,amscd,bbm,amsthm,mathrsfs,dsfont}
\theoremstyle{plain}
\newtheorem{theorem}{Theorem}[section]

\newtheorem{lemma}[theorem]{Lemma}
\newtheorem{prop}[theorem]{Proposition}
\theoremstyle{definition}
\newtheorem{defn}[theorem]{Definition}
\newtheorem{exa}[theorem]{Example}

\theoremstyle{remark}

\newtheorem{remark}{Remark}[section]

\DeclareMathOperator{\Ext}{Ext} 

\begin{document}
\title[Piecewise-Koszul algebras]{\bf Piecewise-Koszul algebras
}
\author{Jia-Feng Lu}
\address{Department of Mathematics, Zhejiang University, Hangzhou
310027, China} \email{jiafenglv@126.com}
\author{ Ji-Wei He}
\address{Institute of Mathematics, Fudan University, Shanghai
200433, China} \email{hejw@fudan.edu.cn}
\author{ Di-Ming Lu}
\address{Department of Mathematics, Zhejiang University, Hangzhou 310027, China}
\email{dmlu@zju.edu.cn} \keywords{($d$-)Koszul algebras,
$A_{\infty}$-algebras, Yoneda-Ext algebras.} \subjclass[2000]{16E05,
16E40, 16S37, 16W50.}
\date{}
\maketitle

\begin{abstract} It is a small step toward the Koszul-type
algebras. The piecewise-Koszul algebras are, in general, a new
class of quadratic algebras but not the classical Koszul ones,
simultaneously they agree with both the classical Koszul and
higher Koszul algebras in special cases. We give a criteria
theorem for a graded algebra $A$ to be piecewise-Koszul in terms
of its Yoneda-Ext algebra $E(A)$, and show an
$A_{\infty}$-structure on $E(A)$. Relations between Koszul
algebras and piecewise-Koszul algebras are discussed. In
particular, our results are related to the third question of
Green-Marcos' \cite{GM}.
\end{abstract}
\vspace{0.8cm}

\section{Introduction}

The Koszul algebra, introduced by S. Priddy in \cite{P}, is one of
quadratic algebras with a linear resolution. Such an algebra may
be understood a positively graded algebra that is ``as close to
semisimple as it can possible be'' (\cite{BGS}). Many nice
homological properties of Koszul algebras have been shown in
research areas of commutative and noncommutative algebras, such as
algebraic topology, algebraic geometry, quantum group, and Lie
algebra (\cite{AG}, \cite{BGS}, \cite{GM1}, etc.). Thirty years
later, motivated by the cubic Artin-Schelter regular algebras,
Berger extended the concept to higher homogeneous algebras
\cite{B}, one can find more discussions under the name $d$-Koszul
algebras in \cite{GMMZ}, or higher Koszul algebras in \cite{HL},
the latter explained Koszulity by $A_\infty$-language.

In this paper, we introduce a new class of Koszul-type algebras, we
name it {\it piecewise-Koszul\/} algebra. Such an algebra is
determined by a pair of parameters $(p, d)$, one shows its
periodicity, and the other one is related to the degree of jump. It
agrees with the classical Koszul algebra when the period equals to
the jumping degree, and goes back to the $d$-Koszul algebra when the
period $p=2$.

What we are interested for these Koszul-type algebras is that in the
case of $d> p\ge 3$. Such algebras provide a new class of quadratic
algebras but not the classical Koszul algebras. The piecewise-Koszul
algebras behave somewhat like the classical Koszul and higher Koszul
algebras: for example, there is a version of the criteria theorem
for a graded algebra to be Koszul in terms of its Yoneda-Ext
algebra. Suppose that $A=\bigoplus_{i\geq 0}A_i$ is generated in
degree zero and one with $A_0$ semisimple, $E(A)=\bigoplus_{i\ge
0}\mbox{Ext}^i_A(A_0, A_0)$ is the Yoneda-Ext algebra of $A$, which
is bigraded by the (\emph{ext, shift})-degree; that is, its $(i,
j)^{th}$ component is $\Ext^i_A(A_0, A_0)_j$. We also call $E(A)$
the Koszul dual of $A$, as other literatures named. The following
results were stated in \cite{GM1}, \cite{HL} and \cite{GMMZ},
respectively:

\vskip2mm
\begin{itemize}
\item  \textsl{$A$ is a Koszul algebra if and only if  ${E(A)}$ is
generated in the ext-degrees $0, 1$ (and $\mbox{Ext}^1_A(A_0,
A_0)=\mbox{Ext}^1_A(A_0, A_0)_{1}$)};
\item  \textsl{$A$ is a $d$-Koszul algebra if and only if  ${E(A)}$ is
generated in  the ext-degrees $0, 1, 2$, and $\mbox{Ext}^2_A(A_0,
A_0)=\mbox{Ext}^2_A(A_0, A_0)_{d}$}.
\end{itemize}

\vskip2mm For our objects, we have

\begin{itemize}
\item  \textsl{$A$ is a piecewise-Koszul algebra if and only if ${E(A)}$
is generated in  the ext-degrees $0, 1, p$, and $\mbox{Ext}^p_A(A_0,
A_0)=\mbox{Ext}^p_A(A_0, A_0)_{d}$} (Theorem \ref{2}).
\end{itemize}
\vskip2mm

The piecewise-Koszul algebras are one of ``$\delta$-Koszul
algebras'' \cite{GM}. As a result, it seems that our objects show a
negative answer to the third question in \cite{GM} provided we may
get a piecewise-Koszul algebra with any given period $p$. These are
discussed in Section 2. We introduce piecewise-Koszul modules in
this section as well.  Denote $\mathcal{PK}(A)$ the category of
piecewise-Koszul modules, we prove that $\mathcal{PK}(A)$ is
extension closed and co-kernels preserved.

Although the piecewise-Koszul algebra $A$ is quadratic in general,
its Yoneda-Ext algebra $E(A)$ admits at least a non-trivial higher
multiplication, which shows another difference between classical
Koszul algebras and piecewise-Koszul algebras, since Keller showed
in \cite{K2} that a quadratic algebra is Koszul if and only if the
higher multiplications on its Koszul dual are trivial. In Section
3, we investigate the $A_{\infty}$-structure on the Koszul dual of
$A$, where $A$ is a piecewise-Koszul algebra. It turns out that
there is a family of nontrivial higher multiplications on $E(A)$.
In particular, we give a criteria for a graded algebra to be a
piecewise-Koszul algebra in terms of considering the Koszul dual
$E(A)$ as a $(2,d\!-\!1)$-algebra.

However, piecewise-Koszul algebras resemble with the classical
Koszul algebras on the other hand. In fact, piecewise-Koszul
algebras are turned out to be closely related with Koszul algebras.
We will discuss their relations in Section 4.

Let $\mathbb{Z}$ and $\mathbb{N}$ denote the set of integers and
natural numbers respectively.

Throughout we work over a fixed field $\mathbb{F}$.

\medskip

\section{Definitions and Properties}

All the graded $\mathbb{F}$-algebras $A=\bigoplus_{i\geq 0}A_i$ are
assumed with the following properties: (a) $A_0$ is a semisimple
Artin algebra; (b) $A$ is generated in degrees 0 and 1; that is,
$A_i \cdot A_j = A_{i+j}$ for all $0\leq i, j<\infty$; (c) $A_1$ is
of finite dimension as an $\mathbb{F}$-space.

The graded Jacobson radical of $A$, denoted by $J$, is
$\bigoplus_{i\geq 1}A_i$. Let $Gr(A)$ denote the category of graded
$A$-modules, and $gr(A)$ its full subcategory of finitely generated
$A$-modules. The morphisms in these categories, denoted by
Hom$_{Gr(A)}(M, N)$, are graded $A$-module maps of degree zero. We
denote $Gr_s(A)$ and $gr_s(A)$ the full subcategories of $Gr(A)$ and
$gr(A)$ whose objects are generated in degree $s$, respectively. An
object in $Gr_s(A)$ or $gr_s(A)$ is called a {\it pure} $A$-module.

Endowed with the Yoneda product, $\hbox{\rm
Ext}^*_A(A_0,A_0)=\bigoplus_{i\ge0}\hbox{\rm Ext}_A^i(A_0,A_0)$ is
an associative graded algebra. Let $M$ and $N$ be finitely
generated graded $A$-modules. Then $\hbox{\rm
Ext}^*_A(M,N)=\bigoplus_{i\ge0}\hbox{\rm Ext}_A^i(M,N)$ is a
graded left $\hbox{\rm Ext}^*_A(N,N)$-module. For simplicity, we
write $E(A)=\hbox{\rm Ext}^*_A(A_0,A_0)$, and
$\mathcal{E}(M)=\hbox{\rm Ext}^*_A(M,A_0)$ that is a graded
$E(A)$-module, we call it the Koszul dual of $M$.

The Koszul dual $\mathcal{E}(M)$ of a finitely generated graded
module $M$ is bigraded; that is, if $[x] \in$ Ext$_A^i(M,A_0)_j $,
we denote the bidegree of $[x]$ as $(i, j)$, call the first degree
\textit {ext}-degree and the second degree \textit {shift}-degree,
respectively. Similarly, the Yoneda-Ext-algebra $E(A)$ of a
positively graded $\mathbb{F}$-algebra is a bigraded algebra.

Given a pair of integers $d$ and $p$ ($d\geq p\ge 2$), we introduce
a function $\delta^d_p: \mathbb{N}\rightarrow\mathbb{N}$ by
$$\delta^d_{p}(n)=\left\{\begin{array}{llll}
\frac{nd}{p},  &  \mbox{if $n \equiv 0 (\textrm{mod} p)$,}\\
\frac{(n-1)d}{p}+1,  &  \mbox{if $n \equiv 1 (\textrm{mod} p)$,}\\
\cdots& \cdots \\
\frac{(n\!-\!p+1)d}{p}+p\!-\!1, \quad & \mbox{if $n\equiv p\!-\!\!1
(\textrm{mod} p)$.}
\end{array}
\right.
$$

\vskip3mm
\begin{defn}
A graded algebra $A=\bigoplus_{i\geq 0}A_i$ is called a
{\it{piecewise-Koszul algebra\/}} if the trivial $A$-module $A_0$
admits a minimal graded projective resolution
$$ \mathbf{P}:\quad  \cdots \rightarrow P_n\rightarrow \cdots \rightarrow
P_1\rightarrow P_0\rightarrow A_0\rightarrow 0,$$ such that each
$P_n$ is generated in degree $\delta^d_p(n)$ for all $n\geq0$.
\end{defn}

The parameter $p$ shows certain periodicity in the resolution above,
and the parameter $d$ is related to a gap between two segments. The
classical Koszul algebras and the higher Koszul algebras are
piecewise-Koszul algebras, which take place of $d=p$ and $p=2$
respectively, we refer to \cite{BGS}-\cite{P} for the details. What
we are interested, in this paper, is the case of $d> p\ge 3$. Such
algebras provide us a new class of quadratic algebras different from
the classical Koszul algebras.

We will give an example of piecewise-Koszul algebra at the end of
this section.

\begin{defn}
Let $A$ be a piecewise-Koszul algebra and $M\in gr(A)$. We call $M$
a {\it{piecewise-Koszul module\/}} if it has a minimal graded
projective resolution of the form
$$\mathbf{Q}:\quad \cdots \rightarrow Q_n \stackrel{f_n}
\rightarrow  \cdots \rightarrow Q_1\stackrel{f_1}\rightarrow
Q_0\stackrel{f_0}\rightarrow  M\rightarrow 0,$$ in which each $Q_n$
is generated in degree $\delta^d_p(n)+s$ for $n\ge 0$, where $s$ is
a fixed integer.
\end{defn}

It is clear, by the definitions above, that a graded
$\mathbb{F}$-algebra $A$ is a piecewise-Koszul algebra (with the
period $p$ and jump degree $d$) if and only if the $n^{th}$
component $\text{\rm Ext}_A^n(A_0,A_0)$ is concentrated in degree
$(n,\delta^d_p(n))$ for all $n\geq 0$. Let $A$ be a piecewise-Koszul
algebra and $M\in gr(A)$, then $M$ is a piecewise-Koszul module if
and only if $\text{\rm Ext}_A^n(M,A_0)$ is concentrated in the
degree $(n, \delta^d_p(n)+s)$ for all $n\geq 0$, where $d, p, s$ are
fixed integers in the definition.

\begin{defn}\cite{LPWZ1} An $A_\infty$-algebra over a field $\mathbb{F}$ is a
$\mathbb{Z}$-graded vector space
$$E=\bigoplus_{n\in \mathbb{Z}}E^n$$endowed with a family of graded
$\mathbb{F}$-linear maps
$$m_n: E^{\otimes n}\longrightarrow E, \;\;(n\geq 1)$$
of degree $(2-n)$ satisfying the following {\it Stasheff's
identities:\/} for all $n\geq 1$,

\bigskip \noindent {\bf SI(n)} $\qquad\qquad \qquad\qquad$ $\sum
(-1)^{r+st} m_u(1^{\otimes r}\otimes m_s \otimes 1^{\otimes t})=0,$

\bigskip
\noindent where the sum runs over all decomposition $n= r+s+t$, $r$,
$t\geq 0$ and $s\geq 1$, and where $u=r+1+t$. Note that when these
formulas are applied to homogeneous elements, additional signs
appear due to the Koszul sign rule. If $E$ is bi-graded, we can also
define bi-graded $A_{\infty}$-algebras similarly and all the
multiplications preserve the second degree, we refer to \cite{LPWZ1}
for the details.

The multiplications $\{m_n\}_{n\geq 3}$ are called the higher
multiplications of $E$.

An {\it{augmented $A_\infty$-algebra\/}} $E$ is an
$A_\infty$-algebra satisfying the following conditions, where for
the notions of {\it{$A_\infty$-morphism\/}} and {\it{strict
$A_\infty$-morphism\/}}, we refer to \cite{HL} and \cite{LPWZ1}
for the details,
\begin{enumerate}
\item There is a strict $A_\infty$-morphism $\eta_E:
\mathbb{F}\rightarrow E$ such that $m_n(1^{\otimes i}\otimes
\eta_E\otimes 1^{\otimes j})=0$ for all $n\neq 2$ and $i+j=n-1$, and
$m_2(1\otimes \eta_E)=m_2(\eta_E\otimes 1)=1_E$; \item There is a
strict $A_\infty$-morphism $\varepsilon_E: E\rightarrow \mathbb{F}$
such that $\varepsilon_E\circ \eta_E=1$.
\end{enumerate}
\end{defn}

 The following definition is related to the Koszul
dual $E(A)$ as a special $A_\infty$-algebras, where $A$ is a
piecewise-Koszul algebra.

\begin{defn}\cite{HL}
Let $E=\bigoplus_{n \in \mathbb{Z}}E^n$ be an $A_{\infty}$-algebra.
If $E$ has only two nontrivial multiplications $m_2$ and $m_l$, then
$(E, m_2, m_l )$ is called a $(2,l)$-algebra.

An augmented $(2,l)$-algebra $(E, m_2, m_l )$ is called a {\it
{reduced $(2,l)$-algebra}\/} provided the following conditions are
satisfied:
\begin{enumerate}
\item $E=\mathbb{F}\oplus E^1\oplus E^2\oplus\cdots$; \item
$m_2(E^{3k_1+s_1}\otimes E^{3k_2+s_2})=0$ for all $3\leq
s_1+s_2\leq 4$ and $s_j= 1$ or $2$ and $k_j\geq 0,\;\; j=1,\;2$;
\item $m_l(E^{i_1}\otimes E^{i_2}\otimes \cdots \otimes
E^{i_l})=0$ unless one of the $i_j\equiv 2(\textrm{mod} 3)$ and the
others $\equiv 1(\textrm{mod} 3)$.
\end{enumerate}
A reduced $(2,l)$-algebra $E$ is said to be generated by $E^1$ if
for all $n\geq 2$,
$$E^n=\sum_{i+j=n; \; i,j\geq 1}m_2(E^i\otimes E^j)+\sum m_l(E^{i_1}\otimes E^{i_2}\otimes \cdots \otimes E^{i_l}),$$
where the sum in the second sigma runs over all decompositions
$i_1+i_2+\cdots +i_l+2-l=n$ and all $i_t\geq 1$ and $t=1, 2, \cdots,
l$.
\end{defn}
In what follows, we use $\Omega^{n}(M)$ to denote the $n^{th}$
syzygy of $M$.

We have the following equivalent descriptions on piecewise-Koszul
objects and Theorem \ref{2} seems to give a negative answer to
Green's open question in \cite{GM}.
\begin{lemma} \cite{GMMZ}\label{2.9} Suppose that $\mathbf{P}$
is the minimal graded projective resolution
 of the trivial $A$-module $A_0$ and $P_n$ is finitely generated
with generators in degree $\delta^d_p(n)$. Assume that
$\delta^d_p(i+j)=\delta^d_p(i)+\delta^d_p(j)$. Then the Yoneda map
$$\text{\rm Ext}_A^i(A_0,A_0)\otimes_{\mathbb{F}}\text{\rm Ext}_A^j(A_0,A_0)\rightarrow
\text{\rm Ext}_A^{i+j}(A_0,A_0)$$ is surjective. Moreover,
$$\Ext^{i+j}_A(A_0,A_0)=\Ext^{i}_A(A_0,A_0)\cdot
\Ext^{j}_A(A_0,A_0)=\Ext^{j}_A(A_0,A_0)\cdot \Ext^{i}_A(A_0,A_0).$$
\hfill{$\square$}\end{lemma}

\begin{theorem}\label{2}
Let $A=\bigoplus_{i\geq 0}A_i$ be a graded $\mathbb{F}$-algebra as
assumed before, and $E(A)$ its Yoneda-Ext-algebra. Then the
following statements are equivalent:

(1) $A$ is a $piecewise$-Koszul algebra;

(2) $E(A)$ is generated in the ext-degrees  $0, 1$ and $p$,
moreover, $\text{\rm Ext}_A^p(A_0,A_0) = \text{\rm
Ext}_A^p(A_0,A_0)_{d}$.
\end{theorem}
\begin{proof} We write the function $\delta^d_p(n)$ as $\delta(n)$ for simplicity.

(1)$\Rightarrow$(2). This follows from Lemma \ref{2.9} and an easy
induction on the ext-degrees.

(2)$\Rightarrow$(1). Suppose that the minimal projective resolution
of $A_0$ is $$ \cdots \rightarrow P_n\rightarrow \cdots \rightarrow
P_1\rightarrow P_0\rightarrow A_0\rightarrow 0.$$

When $0\le n < p$, we have
$\Ext^n_A(A_0,A_0)=\Ext^n_A(A_0,A_0)_{n}=(\Ext^1_A(A_0, A_0))^n$.
Hence $P_n$ is generated in degree $\delta(n)$ for $0\le n < p$.

In general, we shall show that $P_n$ is generated in degree
$\delta(n)$ for all $n\ge0$.

We only prove the case of $k=i=1$. Observing that $\text{\rm
Ext}_A^p(A_0,A_0) = \text{\rm Ext}_A^p(A_0,A_0)_{d}$ and $\text{\rm
Ext}^{s}_A(A_0,A_0)\cdot \text{\rm Ext}^t_A(A_0,A_0)=0$ for $s,t<p$,
we have
\begin{eqnarray*}
  &&\text{\rm
Ext}^{p+1}_A(A_0,A_0) \\&=&\sum_{s+t=p+1;\; s, t \leq p} \text{\rm
Ext}^{s}_A(A_0,A_0)\cdot \text{\rm Ext}^t_A(A_0,A_0)\\&=&
  \sum_{s+t=p+1,\; s, t< p} \text{\rm
Ext}^{s}_A(A_0,A_0)\cdot \text{\rm Ext}^t_A(A_0,A_0)\\&&\;
+\Ext^{1}_A(A_0,A_0)\cdot
\Ext^{p}_A(A_0,A_0)+\Ext^{p}_A(A_0,A_0)\cdot \Ext^{1}_A(A_0,A_0)\\
                             &=& \Ext^{1}_A(A_0,A_0)\cdot
\Ext^{p}_A(A_0,A_0)+\Ext^{p}_A(A_0,A_0)\cdot
\Ext^{1}_A(A_0,A_0)\\
                              &=&  \Ext^{1}_A(A_0,A_0)_1\cdot
\Ext^{p}_A(A_0,A_0)_d+\Ext^{p}_A(A_0,A_0)_d\cdot
\Ext^{1}_A(A_0,A_0)_1\\
&=&\Ext^{p+1}_A(A_0,A_0)_{d+1}.
                          \end{eqnarray*}
Hence $P_{p+1}$ is generated in degree $\delta(p+1)$.
\end{proof}

There is a similar description for piecewise-Koszul modules.
\begin{theorem}\label{4}
Let $A$ be a piecewise-Koszul algebra and $M\in gr_0(A)$. Then $M$
is a piecewise-Koszul module if and only if $\mathcal{E}(M)$ is
generated in degree 0 as a graded $E(A)$-module.
\end{theorem}
\begin{proof} This is an easy consequence of  Proposition 3.5 in
\cite{GMMZ}.
\end{proof}

Let $A$ be a piecewise-Koszul algebra and $M\in gr_s(A)$ be a
piecewise-Koszul module. Then we have the following exact sequences
$$0\rightarrow \Omega^n(M)\rightarrow \Omega^n(M/JM)\rightarrow\Omega^{n-1}(JM)\rightarrow 0,$$
for all $n\geq 1$. All modules in the exact sequences above are
generated in degree $\delta_p^d(n)+s$. We get the following exact
sequences,
$$0\rightarrow \text{\rm Ext}_A^{n-1}(JM,A_0)\rightarrow \text{\rm Ext}_A^{n}(M/JM,A_0)\rightarrow
\text{\rm Ext}_A^{n}(M,A_0)\rightarrow 0.$$ Moreover, in the
shift-grading, all the modules in the exact sequences above are
concentrated in degree $\delta^d_p(n)+s$ for a fixed integer $s$.

The following proposition shows that the category
$\mathcal{PK}(A)$ of piecewise-Koszul modules is closed under
extension and preserves cokernels. We omit the proofs since they
are obvious.
\begin{prop}\label{8}
Let $$0\rightarrow K\rightarrow M\rightarrow N\rightarrow 0$$ be a
short exact sequence in $gr(A)$. Then we have the following
statements,
\begin{enumerate}
\item If $K$ and $N$ are in $\mathcal{PK}(A)$ with $K$ and $N$
being in $gr_s(A)$, then $M\in\mathcal{PK}(A)$,
\item If $K$ and $M$ are in $\mathcal{PK}(A)$ with $K$ and $M$
being in $gr_s(A)$, then $M/K\cong N\in\mathcal{PK}(A)$.
\end{enumerate}\hfill{$\square$}
\end{prop}
\begin{prop}\label{7}
Let $M\in gr(A)$ be a piecewise-Koszul module. Then
\begin{enumerate}
\item All the $(pn)^{th}$ syzygies $\Omega^{pn}(M)$ are
piecewise-Koszul modules;
\item All the $(pn-1)^{th}$ syzygies $\Omega^{pn-1}(JM)$ are
piecewise-Koszul modules.
\end{enumerate}
\end{prop}
\begin{proof}
The first assertion is immediate from the definition of the
piecewise-Koszul modules. For the second statement, clearly, we have
the following exact sequences in $gr_{\delta^d_p(pn)+s}(A)$
$$0\rightarrow \Omega^{pn}(M)\rightarrow
\Omega^{pn}(M/JM)\rightarrow \Omega^{pn-1}(JM)\rightarrow 0,$$ by
assertion (1), $\Omega^{pn}(M)$ and $\Omega^{pn}(M/JM)$ are in
$\mathcal{PK}(A)$ since $M$ and $M/JM$ are in $\mathcal{PK}(A)$. By
Proposition \ref{8}, we have that $\Omega^{pn-1}(JM)\in
\mathcal{PK}(A)$.
\end{proof}
\begin{exa}
Let $\mathbb{F}$ be a field and let $\Gamma$ be the quiver:

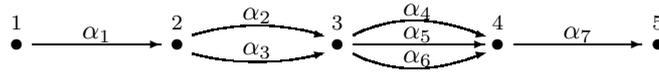
\unitlength=1.2pt
\begin{figure}[h]
\begin{picture}(200,28)(0,0)
\put(0,10){\circle*{3}}\put(0,17){\makebox(0,0){\footnotesize$1$}}
\put(5,10){\vector(1,0){40}}\put(25,13){\makebox(0,0){$\alpha_1$}}
\put(50,10){\circle*{3}}\put(50,17){\makebox(0,0){\footnotesize$2$}}
\put(75,19){\makebox(0,0){$\alpha_2$}}
\put(96.2,12.7){$\vector(4,-1){0}$}
\put(75,8){\makebox(0,0){$\alpha_3$}}
\put(96.2,7.5){$\vector(4,1){0}$}
\put(100,10){\circle*{3}}\put(100,17){\makebox(0,0){\footnotesize$3$}}
\put(125,20.2){\makebox(0,0){$\alpha_4$}}
\put(146.5,12.7){$\vector(4,-1){0}$}
\put(105,10){\vector(1,0){41}}\put(125,13){\makebox(0,0){$\alpha_5$}}
\put(125,5.6){\makebox(0,0){$\alpha_6$}}
\put(146.5,7.4){$\vector(4,1){0}$}
\put(150,10){\circle*{3}}\put(150,17){\makebox(0,0){\footnotesize$4$}}
\put(155,10){\vector(1,0){40}}\put(175,13){\makebox(0,0){$\alpha_7$}}
\put(200,10){\circle*{3}}\put(200,17){\makebox(0,0){\footnotesize$5$}}
\qbezier(55,13)(75,19)(95,13)\qbezier(55,7)(75,2)(95,7)
\qbezier(105,7)(123,-2)(145,7)\qbezier(105,13)(123,22)(145,13)
\end{picture}
\end{figure}

Now let $A={\mathbb{F}\Gamma}/{R}$, where $R$ is the ideal generated
by the following relations: $$ \alpha_1\alpha_2-\alpha_1\alpha_3
,\quad \alpha_4\alpha_7-\alpha_5\alpha_7,\quad
\alpha_5\alpha_7-\alpha_6\alpha_7,\quad \alpha_2\alpha_4,\quad
\alpha_3\alpha_6.$$ It is not difficult to check that $A$ is a
piecewise-Koszul algebra with $p=3$ and $d=4$.
\end{exa}

\medskip
\section{The $A_{\infty}$-structure on the Koszul dual $E(A)$}
Let $A$ be a piecewise-Koszul algebra and $E(A)$ its Koszul dual. We
will investigate the $A_\infty$-structure on $E(A)$ in detail in
this section. As a result, a criteria for a graded algebra to be
piecewise-Koszul in terms of the $A_{\infty}$-algebra is given.

For the simplicity, we assume, in this section, that
$A_0=\mathbb{F}$ in the graded algebra $A$ and the piecewise-Koszul
objects are related to $\delta^d_3(n)$ with $d> 3$ a fixed integer.

\begin{lemma}(\cite{K2}, \cite{HL})\label{5.1}
Let $A$ be an augmented differential bi-graded algebra and $E=HA$.
Then there is an augmented $A_{\infty}$-structure $\{m_i\}$ on $E$
such that $m_1=0$, $m_2$ is induced by the multiplication of $A$,
and $E$ is quasi-isomorphic to $A$ as augmented bi-graded
$A_{\infty}$-algebras.
\end{lemma}

Let $A=\mathbb{F}\oplus A_1\oplus A_2\oplus \cdots$ be a graded
algebra with $dim_{\mathbb{F}}(A_i)< \infty$ for all $i$. Let $J$
be the graded Jacobson radical of $A$, $J^*$ be the graded dual of
$J$ and $T(J^*)$ the tensor algebra. From (\cite{K2}, \cite{HL}),
we know that $T(J^*)$ is an augmented differential bi-graded
algebra and $E(A)=\Ext_A^*(\mathbb{F},\mathbb{F})\cong H(T(J^*))$.
By Lemma \ref{5.1}, the Koszul dual of $A$, $E(A)$, is an
augmented bi-graded $A_{\infty}$-algebra up to quasi-isomorphism.
We will call such $A_{\infty}$-structure on $E(A)$ induced from
$T(J^*)$.

\begin{prop}\label{5.2}
Let $A$ be a piecewise-Koszul algebra and $E(A)$ be its Koszul dual.
Then
$\Ext^3_A(\mathbb{F},\mathbb{F})=\Ext^3_A(\mathbb{F},\mathbb{F})_d$,
and all possible $A_\infty$-structures $\{m_q\}$ on
 $E(A)$ satisfy $m_p=0$ if $q\neq k(d-3)+2$ for all
$k\in \mathbb{N}$. Moreover, $(E(A), \{m_q\})$ satisfies the
conditions (1) and (2) in the Definition 2.4.
\end{prop}
\begin{proof}
Since $A$ is a piecewise-Koszul algebra, we have
$E^{3k+s}:=\Ext_A^{3k+s}(\mathbb{F},\mathbb{F})=E^{3k+s}_{kd+s}$
where $s=0,\ 1$ or $2$. Recall that all the multiplications
$\{m_q\}$ preserve the second degree, we have
$$m_q(E^{3k_1+s_1}\otimes E^{3k_2+s_2}\otimes\cdots \otimes E^{3k_q+s_q})\subset
E^{3(k_1+k_2+\cdots k_q)+s_1+s_2+\cdots +s_q+2-q}_{(k_1+k_2+\cdots
k_q)d+s_1+s_2+\cdots +s_q},$$where $s_i=0$, $1$ or $2$ $(1\leq i\leq
q)$. When $s_1+s_2+\cdots +s_q+2-q=3k$ and if $m_q\neq 0$, then
$s_1+s_2+\cdots +s_q=kd$ and $q=k(d-3)+2$. For the cases
$s_1+s_2+\cdots +s_q+2-q=3k+1$ and $s_1+s_2+\cdots +s_q+2-q=3k+2$,
similarly, we can get the same result. For the last statement, it is
clear that $E^0=\mathbb{F}$. Notice that $d>3$, we have
$$m_2(E^{3k_1+1}_{k_1d+1}\otimes E^{3k_2+2}_{k_2d+2})\subset E^{3(k_1+k_2)+3}_{(k_1+k_2)d+3}=0,$$
$$m_2(E^{3k_1+2}_{k_1d+2}\otimes E^{3k_2+1}_{k_2d+1})\subset
E^{3(k_1+k_2)+3}_{(k_1+k_2)d+3}=0,$$and
$$m_2(E^{3k_1+2}_{k_1d+2}\otimes E^{3k_2+2}_{k_2d+2})\subset E^{3(k_1+k_2)+4}_{(k_1+k_2)d+4}=0.$$
\end{proof}

If $E$ is an $A_{\infty}$-algebra with nontrivial multiplications
$m_2$ and $m_l$, then it is easy to see that all the Stasheff
identities hold automatically except for the following three cases:
\begin{equation}  m_2(m_2\otimes 1)=m_2(1\otimes m_2),
\tag*{$\mathbf{SI(3)}$:}\end{equation}
\begin{equation} \sum_{i+j=l-1}(-1)^{i+lj}m_l(1^{\otimes i}\otimes m_l\otimes
1^{\otimes j})=0,\tag*{$\mathbf{SI(2l-1)}$:}\end{equation} \noindent
and
\begin{equation}
\sum_{i+j=l-1}(-1)^{i+lj}m_l(1^{\otimes i}\otimes m_2\otimes
1^{\otimes j})=m_2(1\otimes m_l)-(-1)^{l}m_2(m_l\otimes
1).\tag*{$\mathbf{SI(l+1)}$:}\end{equation}

From \cite{K1}, we have the following lemma.
\begin{lemma}\label{5.5}
Let $A=\mathbb{F}\oplus A_1\oplus A_2\oplus \cdots$ be a graded
algebra generated in degree 1. Then there exits an $A_{\infty}$-
structure on $\Ext_A^*(\mathbb{F},\mathbb{F})$, such that
$\Ext_A^*(\mathbb{F},\mathbb{F})$ is generated by
$\Ext_A^1(\mathbb{F},\mathbb{F})$ as an $A_{\infty}$-algebra.
\end{lemma}
\begin{remark}
By Lemma \ref{5.5}, it is easy to see that there exists at least
one nontrivial higher multiplication constructed in Proposition
\ref{5.2}.
\end{remark}
\begin{theorem}\label{5.6}
Let $A=\mathbb{F}\oplus A_1\oplus A_2\oplus \cdots$ be a graded
algebra generated in degree 1 and $E(A)$ be its Koszul dual.
\begin{enumerate}
\item If $A$ is a piecewise-Koszul algebra, then there exists an
$A_{\infty}$-structure $\{m_q\}$ (the same as Proposition \ref{5.2})
such that $(E(A), \{m_q\})$ is an $A_{\infty}$-algebra generated by
$E^1(A)$; \item Conversely, if the nontrivial multiplications in the
$A_{\infty}$-structure on $E(A)$ induced from $T(J^*)$ are $m_2$ and
$m_{d-1}$, such that $(E(A),m_2,m_{d-1})$ is a reduced
$(2,d-1)$-algebra generated by $E^1(A)$. Then $A$ is a
piecewise-Koszul algebra with jumping degree $d$.
\end{enumerate}
\end{theorem}
\begin{proof}
The statement (1) follows from Proposition \ref{5.2} and Lemma
\ref{5.5}. Since $A$ is generated by $A_1$, we have
$E^1(A)=E^1(A)_1$. Therefore $E^2(A)=E^2(A)_2$, since
$E^2(A)=m_2(E^1(A)\otimes E^1(A))$. Also
$$E^3(A)=\sum_{i_1+i_2+\cdots+i_{d-1}+2-(d-1)=3,\;\;i_s\geq 1}
m_{d-1}(E^{i_1}(A)\otimes E^{i_2}(A)\otimes \cdots \otimes
E^{i_{d-1}}(A))$$  implies that $E^3(A)=E^3(A)_{\delta_3^d(3)}$.
Now assume that we have the identities:
$E^n(A)=E^n(A)_{\delta_3^d(n)}$ for all $n< 3k$, where $k\in
\mathbb{N}$. Now consider $E^{3k}(A)$, $E^{3k+1}(A)$ and
$E^{3k+2}(A)$ respectively. By hypothesis,
\begin{eqnarray*}
&&E^{3k}(A)=\sum_{i+j=3k; \; i,j\geq 1}m_2(E^i(A)\otimes
E^j(A))\\&+&\sum_{i_1+i_2+\cdots+i_{d-1}+2-(d-1)=3k,\;\;i_s\geq 1}
m_{d-1} (E^{i_1}(A)\otimes E^{i_2}(A)\otimes \cdots \otimes
E^{i_{d-1}}(A))\\&=&\sum_{3p+3q=3k; \; p,\;q\geq
1}m_2(E^{3p}(A)\otimes
E^{3q}(A))\\&+&\sum_{i_1+i_2+\cdots+i_{d-1}+2-(d-1)=3k,\;\;i_s\geq
1} m_{d-1} (E^{i_1}(A)\otimes E^{i_2}(A)\otimes \cdots \otimes
E^{i_{d-1}}(A)),\\&=& E^{3k}(A)_{\delta_3^d(3k)}.
\end{eqnarray*}
\begin{eqnarray*}
&&E^{3k+1}(A)=\sum_{i+j=3k+1; \; i,j\geq 1}m_2(E^i(A)\otimes
E^j(A))\\&+&\sum_{i_1+i_2+\cdots+i_{d-1}+2-(d-1)=3k+1,\;\;i_s\geq
1} m_{d-1} (E^{i_1}(A)\otimes E^{i_2}(A)\otimes \cdots \otimes
E^{i_{d-1}}(A))\\&=&\sum_{p+q=k; \; p,\;q\geq
1}m_2(E^{3p+1}(A)\otimes E^{3q}(A)+E^{3p}(A)\otimes
E^{3q+1}(A))\\&=& E^{3k+1}(A)_{\delta_3^d(3k+1)}.
\end{eqnarray*}
And similarly, we can prove
$E^{3k+2}(A)=E^{3k+2}(A)_{\delta_3^d(3k+2)}$. Therefore, $A$ is a
piecewise-Koszul algebra with jumping degree $d$.
\end{proof}

\section{Koszul objects from piecewise-Koszul objects}
 We will construct Koszul objects from
the given piecewise-Koszul objects in this section.

Let $A$ be a piecewise-Koszul algebra and $M\in gr_0(A)$ be a
piecewise-Koszul module respect to $\delta_3^d(n)$. For  $k\geq
1$, we set
$$\mathbf{E_k}(A)=\bigoplus_{n\geq 0}\text{\rm Ext}_A^{3kn}(A_0,A_0)$$ and
$$\mathbf{E_k}(M)=\bigoplus_{n\geq 0}\text{\rm Ext}_A^{3kn}(M,A_0).$$
Naturally $\mathbf{E_k}(A)$ is a subalgebra of $E(A)$, and
$\mathbf{E_k}(M)$ is a graded $\mathbf{E_k}(A)$-module.
\begin{theorem}\label{2.8}
Let $A$ be a piecewise-Koszul algebra and $M$ be a
piecewise-Koszul module respect to $\delta_3^d(n)$. The notations
$\mathbf{E_k}(A)$ and $\mathbf{E_k}(M)$ are as defined above. Then
for all integers $k\geq 1$, we have
\begin{enumerate}
\item $\mathbf{E_k}(M)= \bigoplus_{n\geq
0}\Ext_A^{3kn}(M,A_0)$ is a Koszul $\mathbf{E_k}(A)$-module, and
\item $\mathbf{E_k}(A)= \bigoplus_{n\geq
0}\Ext_A^{3kn}(A_0,A_0)$ is a Koszul algebra.
\end{enumerate}
\end{theorem}
\begin{proof}
By Theorem \ref{4}, we get that $\mathbf{E_k}(M)$ is generated in
degree 0 as a graded $\mathbf{E_k}(A)$-module.

Clearly, we have the following exact sequences, for all $n,k>0$,
$$0\rightarrow \text{\rm Ext}_A^{3kn-1}(JM,A_0)\rightarrow\text{\rm Ext}_A^{3kn}(M/JM,A_0)\rightarrow
\text{\rm Ext}_A^{3kn}(M,A_0)\rightarrow 0$$ such that all the
modules in the above exact sequences are concentrated in degree
$\delta_3^d(3nk)$ in the shift-grading.

We have the following exact sequences
$$0\rightarrow \text{\rm Ext}_A^{3k(n-1)}(\Omega
^{3k-1}(JM),A_0)\rightarrow\text{\rm
Ext}_A^{3kn}(M/JM,A_0)\rightarrow \text{\rm
Ext}_A^{3kn}(M,A_0)\rightarrow 0$$ since
$$\text{\rm Ext}_A^{3kn-1}(JM,A_0)=\text{\rm Ext}_A^{3k(n-1)}(\Omega
^{3k-1}(JM),A_0).$$ By taking the direct sums of  the above exact
sequences, we have {\small{$$\begin{array}{c}
                                0\rightarrow \mathbf{E_k}(\Omega
^{3k-1}(JM))[1]
                                \rightarrow
\bigoplus_{n> 0}\text{\rm
Ext}_A^{3kn}(M/JM,A_0)\rightarrow\bigoplus_{n> 0} \text{\rm
Ext}_A^{3kn}(M,A_0)\rightarrow 0.
                              \end{array}$$}}

Now we claim that $\mathbf{E_k}(M/JM)$ is a projective cover of
$\mathbf{E_k}(M)$ and it is generated in degree 0. In fact,
$\mathbf{E_k}(M/JM)$ is a $\mathbf{E_k}(A)$-projective module since
$M/JM$ is semi-simple. $M/JM$ is a piecewise-Koszul module since $A$
is a piecewise-Koszul algebra. Hence $\mathbf{E_k}(M/JM)$ is
generated in degree 0 as a graded $\mathbf{E_k}(A)$-module, and by
the above exact sequence, it is the graded projective cover of
$\mathbf{E_k}(M)$.

Therefore the first syzygy is $\bigoplus_{n> 0}\text{\rm
Ext}_A^{3k(n-1)}(\Omega ^{3k-1}(JM),A_0)$, and $\Omega
^{3k-1}(JM)$ is generated in degree $\delta^d_3(3k)$ and by
Proposition 2.9, $\Omega ^{3k-1}(JM)$ is again a piecewise-Koszul
module. Inductively, we prove the first assertion.

The second assertion is the case of $M=A_0$.  The proof is finished.
\end{proof}
\begin{remark}
If the piecewise-Koszul objects are respect to $\delta_p^d(n)+s$ in
the above theorem, then the theorem can be restated as follows: Let
$A$ be a piecewise-Koszul algebra and $M$ be a piecewise-Koszul
module. Then for all integers $k\geq 1$, we have
\begin{enumerate}
\item $\mathbf{\mathcal{E}_k}(M)= \bigoplus_{n\geq 0}\text{Ext}_A^{pkn}(M,A_0)$ is a
Koszul module, and
\item $\mathbf{\mathcal{E}_k}(A)= \bigoplus_{n\geq 0}\text{Ext}_A^{pkn}(A_0,A_0)$
is a Koszul algebra.
\end{enumerate}
\end{remark}

By Theorem \ref{2.8}, we can construct new Koszul objects from the
given piecewise-Koszul objects. However, we can't construct
$d$-Koszul objects from a given piecewise-Koszul object for $d> 2$.

\vskip3mm

\bigskip
{\it Acknowledgment} $\quad$ Di-Ming Lu is supported by the NSFC
(project 10571152).

\vspace{5mm}

\bibliographystyle{amsplain}

\begin{thebibliography}{11}

\bibitem{AG} R. M. Aquino and E. L. Green, \emph{On modules with linear
presentations over Koszul algebras}, Comm. Algebra \textbf{33}
(2005) 19-36.
\bibitem{ARS} M. Auslander, I. Reiten and S. Smal{\o},
\emph{Representation Theory of Artin Algebras}, Cambridge studies in
Advanced Mathethematics, Vol. \textbf{36}, Cambridge. UK: Cambridge
University Press (1995).
\bibitem{AS} M. Artin and W. F. Schelter, \emph{Graded algebras of global
dimension 3}, Adv. math., Vol. \textbf{66} (1987), 171-216.

\bibitem{BGS} A. Beilinson, V. Ginszburg,
and W. Soergel, \emph{Koszul duality patterns in representation
theory}, J. Amer. Math. Soc., Vol. \textbf{9} (1996), 473-525.
\bibitem{B} R. Berger,
\emph{Koszulity for nonquadratic algebras}, J. Alg., Vol.
\textbf{239} (2001), 705-734.
\bibitem{GM} E. L. Green, E. N. Marcos, \emph{$\delta$-Koszul algebras},
Comm. Alg., Vol. \textbf{33}(6) (2005), 1753-1764.
\bibitem{GM1} E. L. Green, R. Martinez-Villa, \emph{Koszul and Yoneda
algebras}, Representation theory of algebras (Cocoyoc, 1994), CMS
Conference Proceedings, Vol. \textbf{18}, American Mathematical
Society, Providence, RI, (1996), 247-297.
\bibitem{GMMZ} E. L. Green, E. N. Marcos, R. Martinez-Villa, Pu Zhang, \emph{$D$-Koszul algebras}, J. pure
and Appl. Algebra \textbf{193} (2004), 141-162.
\bibitem{GMR} E. L. Green, R. Martinez-Villa, I. Reiten, $\phi$. Solberg, D. Zacharia, \emph{On
modules with linear presentations}, J. Alg., \textbf{205}(2) (1998),
578-604.
\bibitem{HL} J.-W. He, D.-M. Lu, \emph{Higher Koszul Algebras and
A-infinity Algebras}, J. Alg., \textbf{293} (2005), 335--362.
\bibitem{K1} B. Keller, \emph{A-infinity algebras in representation
theory}, Contribution to the proceedings of ICRA IX, Beijing 2000.
\bibitem{K2} B. Keller, \emph{Introduction to A-infinity algebras
and modules}, Homology Homotopy appl., \textbf{3} (2001),
(electronic), 1-35.
\bibitem{LPWZ1} D.-M Lu, J. H. Palmieri, Q.-s. Wu and J. J.Zhang,
\emph{$A_{\infty}$-algebras for ring theorists}, Alg. Colloq., Vol.
11(1), 91-128, 2004.
\bibitem{P} S. Priddy, \emph{Koszul resolutions},
Trans. Amer. Math. Soc., Vol. \textbf{152} (1970), 39-60.
\bibitem{W} C. A. Weible, \emph{An Introduction to Homological Algebra}, Cambridge Studies in
Avanced Mathematics, Vol. \textbf{38}, Cambridge Univ. Press,
(1995).
\end{thebibliography}

\end{document}